\documentclass[A4,12pt,conference]{smart_2025_Paper}

\usepackage{times}

\usepackage{graphicx}
\usepackage{amsmath}
\usepackage{amsfonts}
\usepackage{amssymb}
\usepackage[numbers]{natbib}
\usepackage{hyperref}
\usepackage[detect-all]{siunitx}
\usepackage{booktabs}
\usepackage{float}
\usepackage{hyperref}

\usepackage{tikz}
\usepackage{xparse} 

\DeclareMathOperator*{\argmin}{arg\,min}

\newcommand{\ith}[1]{^{(#1)}}
\newcommand{\one}{\ith{1}}
\newcommand{\two}{\ith{2}}
\newcommand{\thetai}{\theta\ith{i}}

\newcommand{\ex}{\boldsymbol e_x}
\newcommand{\ey}{\boldsymbol e_y}
\newcommand{\ez}{\boldsymbol e_z}

\newcommand{\xv}{\boldsymbol x}
\newcommand{\Xv}{\boldsymbol X}
\newcommand{\uv}{\boldsymbol u}
\newcommand{\slope}{\boldsymbol s}

\newcommand{\body}{\mathcal B}

\newcommand{\rnd}{\bar r}
\newcommand{\cnd}{\bar c}
\newcommand{\knd}{\bar k}
\newcommand{\Fref}{F_{\rm ref}}
\newcommand{\lref}{l_{\rm ref}}

\newcommand{\D}{X}
\newcommand{\Y}{Y}
\newcommand{\reals}{\mathbb R}

\newcommand{\ideal}{y_{\rm ideal}}
\newcommand{\ideali}{y_{{\rm ideal},i}}
\newcommand{\idealj}{y_{{\rm ideal},j}}

\newcommand{\nadir}{y_{\rm nadir}}
\newcommand{\nadiri}{y_{{\rm nadir},i}}
\newcommand{\nadirj}{y_{{\rm nadir},j}}

\newcommand{\weight}{w}
\newcommand{\wpseudo}{\weight_{\rm pseudo}}
\newcommand{\wpseudoi}{\weight_{{\rm pseudo},i}}

\NewDocumentCommand{\drawSmoothThickCurve}{m m O{gray!70}}{%
  \begin{scope}
    \draw[white, line width=2pt]
      plot[smooth] coordinates {#1}
      --
      plot[smooth] coordinates {#2}
      -- cycle;
    \fill[#3]
      plot[smooth] coordinates {#1}
      --
      plot[smooth] coordinates {#2}
      -- cycle;
  \end{scope}
}

\NewDocumentCommand{\drawSmoothThickCurveShadow}{m m O{gray!30}}{%
  \begin{scope}
    \fill[#3]
      plot[smooth] coordinates {#1}
      --
      plot[smooth] coordinates {#2}
      -- cycle;
  \end{scope}
}
\graphicspath{{./figures/}}

\title{Efficient Design of Compliant Mechanisms Using Multi-Objective Optimization}

\author{ALEXANDER HUMER$^{*}$ AND SEBASTIAN PLATZER$^{*}$}

\heading{Alexander Humer and Sebastian Platzer}

\address{
$^{*}$ Institute of Technical Mechanics, Johannes Kepler University Linz\\
Altenberger Str. 69, A-4040 Linz, Austria\\
e-mail: alexander.humer@jku.at - Web page: www.jku.at/tmech}

\abstract{
    Compliant mechanisms achieve motion through elastic deformation. In this work, we address the synthesis of a compliant cross-hinge mechanism capable of large angular strokes while approximating the behavior of an ideal revolute joint. To capture the competing demands of kinematic fidelity, rotational stiffness, and resistance to parasitic motion, we formulate a multi-objective optimization problem based on kinetostatic performance measures. A hybrid design strategy is employed: an efficient beam-based structural model enables extensive exploration of a high-dimensional design space using evolutionary algorithms, followed by fine-tuning with high-fidelity three-dimensional finite element analysis. The resulting Pareto-optimal designs reveal diverse geometric configurations and performance trade-offs.
}


\begin{document}

\section{\MakeUppercase{Introduction}}

\emph{``If something bends to do what it is meant to do, then it is compliant. If the flexibility that allows it to bend also helps it to accomplish something useful, then it is a compliant mechanism''}, is how Larry Howell describes the essence of compliant mechanisms~\cite{howell2001}. 
Their functionality relies on flexible deformation, which is determined by geometric shapes and material properties. 


A \emph{hinge} can arguably be regarded as the quintessential mechanism—embodying the concept of a revolute joint in kinematic terms. 
Replicating a hinge's functionality using compliant mechanisms first requires identifying its defining characteristics.
These characteristics will, in turn, inform the objectives for the synthesis process.
At its core, the purpose of a hinge is to allow two bodies to rotate relative to one another about a single common axis.
In the ideal case of a perfect hinge, this axis remains fixed with respect to both bodies, and the hinge offers no resistance to rotation, regardless of the amplitude. Conversely, it provides infinite stiffness to all other modes of deformation. These contrasting requirements—flexibility in one direction and rigidity in others—are what makes a hinge a good example for the synthesis of compliant mechanisms based on multi-objective optimization.

\begin{figure}[h!]
    \centering
    \includegraphics[height=3cm]{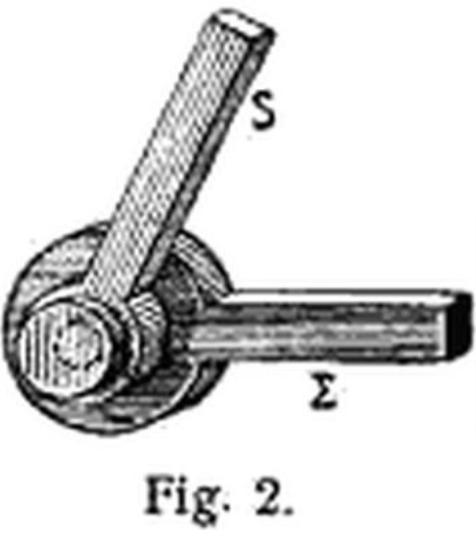}
    \includegraphics[height=3cm]{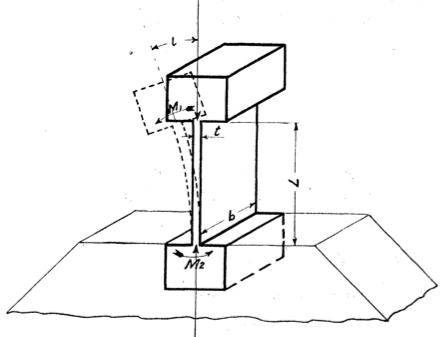}
    \includegraphics[height=3cm]{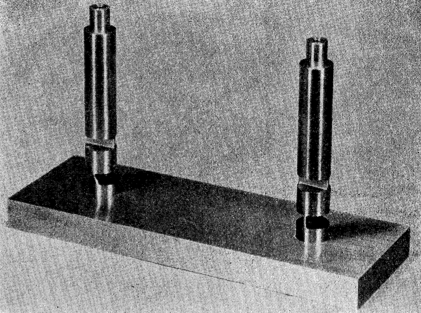}
    \includegraphics[height=3cm]{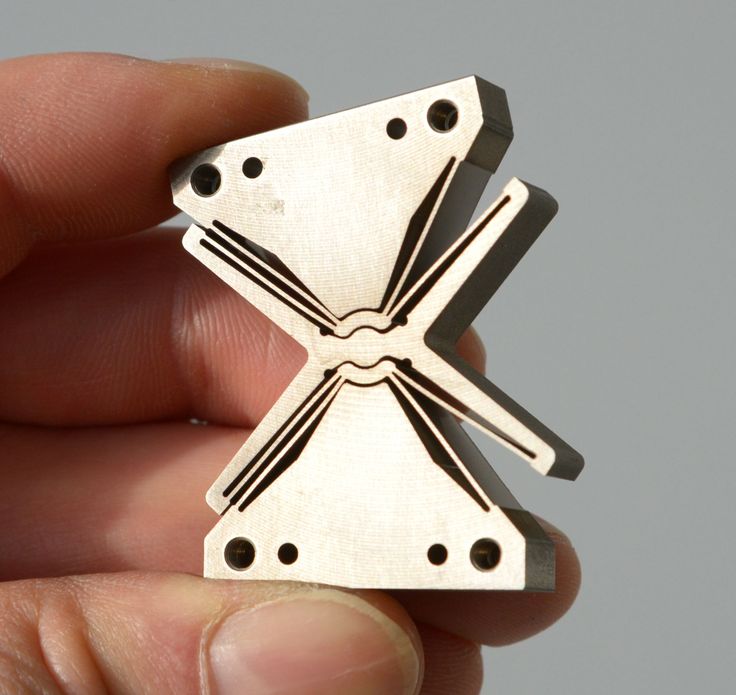}
    \\
    \begin{minipage}[t][1.5em][t]{0.2\textwidth} \small \centering (a) \end{minipage}
	\begin{minipage}[t][1.5em][t]{0.25\textwidth} \small \centering (b) \end{minipage}
	\begin{minipage}[t][1.5em][t]{0.225\textwidth} \small \centering (c) \end{minipage}
	\begin{minipage}[t][1.5em][t]{0.225\textwidth} \small \centering (d) \end{minipage}
    \vspace{-1ex}
    \caption{
        Different kinds of hinges: Conventional pin-hole realization and compliant mechanisms.
    }
    \label{fig:hinges}
\end{figure}

Conventional realizations of revolute joints can be remarkably simple.
Consider, for instance, the pairing of a pin and hole as illustrated in Fig.~\ref{fig:hinges} (a) taken from the encyclopedia~\cite{lueger1904}. 
Relative rotation is enabled by the frictional sliding between the cylindrical surfaces of the pin and hole. 
The most fundamental way to join two (ideally rigid) bodies via a compliant element is to connect them using a thin flexible strip, i.e., a \emph{'flexure'}, whose ends are fixed to each link it connects, see Fig.\ref{fig:hinges} (b).
The idea of using flexures as connecting elements in machinery is by no means a recent concept. 
In arguably one of the first articles on design rules for compliant mechanisms, Eastman~\cite{eastman1937} mentioned wind-tunnel balances of the 19th century using flexure pivots.
None of the variants of flexure hinges illustrated in Fig.\ref{fig:hinges} succeeds in reproducing the behavior of a revolute joint from the kinematic point of view.
A single slender flexure (b) allows for comparatively large rotations among the connected links. 
However, the instantaneous axis of rotation shifts significantly as the deformation—or \emph{stroke}—increases. 
Short flexures (c) exhibit less shift as far as the axis of rotation is concerned. The strength of materials, however, limits their admissible stroke as the flexure’s (effective) slenderness decreases.
To improve kinematics and increase stiffness, multiple flexures can be integrated within a mechanism as, e.g., in CSEM's `butterfly-hinge'~\cite{henein2003} (d). 
Its stroke, however, is limited to \SI{\pm15}{\degree}. 

In a recent survey, Ling et al.~\cite{ling2020} have given a comprehensive overview of three decades of research on the kinetostatic and dynamic modeling of compliant mechanisms. 
The progress in computational methods is reflected in the methodologies, which evolved from engineering beam theories to 3D-finite element analysis.
Optimization lies at the heart of any synthesis problem, irrespective of whether conventional rigid-body or compliant mechanisms are concerned. 
In compliant mechanisms, however, we deal with multiple—typically conflicting—kinematic, static, and dynamic objectives that do not admit globally optimal solutions but demand for compromises. 
Finite stiffness is a defining property of compliant mechanisms, which demands for a trade-off between a maximum compliance regarding the intended motion and minimal compliance in all orthogonal directions. 
By employing \emph{multi-objective optimization} methods and the pivotal notion of \emph{Pareto optimality}, we avoid an a priori weighting among objectives, which inevitably introduces bias for certain designs. 
In what follows, we propose a hybrid two-stage approach to efficiently synthesize mechanisms composed of slender flexures: Structural models provide effective means to capture the global behavior of mechanisms in a first step. 
In a second step, continuum models allow us to resolve details and fine-tune design variants.

\section{CROSS-HINGE GEOMETRY} \label{sec:geometry}
Our focus lies on mechanisms with \emph{distributed compliance} since we aim at large strokes. 
To this end, we generalize the well-established two-flexure cross-hinge, a classical design within the field of compliant mechanisms.
As opposed to the classical cross-hinge, however, we allow for curved flexures with a preferably general shape, see Fig.~\ref{fig:geometry_2d}.
In the undeformed configuration, the cross-sections' angle of rotation $\thetai$ follows a cubic polynomial in a normalized arc-length coordinate $s$,
\begin{equation}
    \thetai (s) = \left(1 - s\right) \thetai_0 + s \, \thetai_1 + 4 s \left(1 - s\right) \thetai_2 + 4 s \left(1 - s\right) \left(2s - 1\right) \thetai_3 ,
\end{equation}
where superscripts $i=1,2$ identify the respective flexure of the cross-hinge.

We use hierarchical polynomials inspired by finite-element shape functions, which guarantees that coefficients $\thetai_0$ and $\thetai_1$ directly correspond to the angles at the terminal points.
We assume cross-sections to be initially perpendicular to the flexures' centerlines.
The position of the centerline $\Xv\ith{i}$ with respect to some global Cartesian frame, whose origin coincides with the point $s=0$ of the first flexure, i.e., $\Xv\ith{1}(0) = \mathbf 0$, then follows upon integration of the unit tangent $\slope\ith{i}$, which is scaled by the length $l\ith{i}$ of the respective flexure, 
\begin{equation}
    \Xv\ith{i} (s) = \Xv\ith{i}_0 + l\ith{i} \int_0^s \slope\ith{i} (\sigma) \mathrm d\sigma, \qquad \slope\ith{i} (s) = \ex \cos \thetai (s) + \ey \sin \thetai (s) ,
\end{equation}
where $\Xv\ith{i}_0$ denotes the position at $s=0$.
For the sake of simplicity, we assume the cross-sections of the flexures to have rectangular shapes, which remain constant along the centerlines;
their widths and heights are denoted by $w\ith{i}$ and $h\ith{i}$, respectively.

\begin{minipage}[c]{0.95\linewidth}
    \centering
    \begin{minipage}[c]{0.45\linewidth}
        \begin{figure}[H]
            \centering
            \includegraphics[scale=0.8]{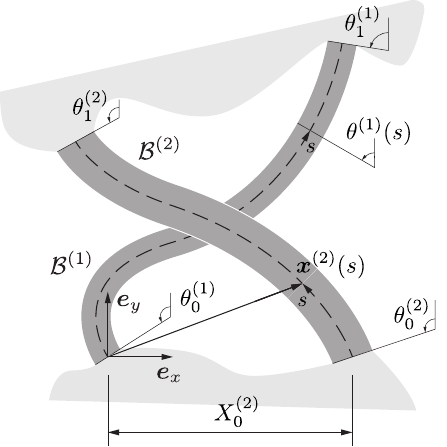}
            \caption{Schematic representation and kinematic description of the shape of the cross-hinge.}
            \label{fig:geometry_2d}
        \end{figure} 
    \end{minipage}
    \hspace{1em}
    \begin{minipage}[c]{0.5\linewidth}
        \begin{table}[H]
            \centering
            \small
            \sisetup{range-phrase=~\ldots~}
            \caption{Admissible ranges of the design variables.}
            \begin{tabular}{r|r||r|r}
                \toprule
                variable	&	\multicolumn{1}{c||}{range}   &   variable &  \multicolumn{1}{c}{range}  \\
                \midrule
                $\theta_0\ith{i} / 1$ 	&	$0 \ldots \pi$ 	&	$\alpha / 1$	&	\numrange{0.5}{2}\phantom{0} \\
                $\theta_1\ith{i} / 1$	&	$-\pi \ldots \pi$ 	&	$\beta\ith{i} / 1$ 	&	\numrange{5}{20} \\
                $\theta_2\ith{i} / 1$	&	$-\pi \ldots \pi$ 	&	$\gamma / 1$ 	&	\numrange{0.5}{2}\phantom{0} \\
                $\theta_3\ith{i} / 1$	&	$-\pi \ldots \pi$ 	&	$\delta / 1$	&	\numrange{0}{1}\phantom{0} \\
                \bottomrule			
            \end{tabular}
            \label{tab:range-variables}
        \end{table}
    \end{minipage}
\end{minipage}

%

The choice of appropriate design variables plays a crucial role in any synthesis problem, and compliant mechanisms are no exception. 
To obtain preferably general results, we introduce a set of non-dimensional similarity parameters as design variables that govern the geometry of the generalized cross-hinge.
The coefficients $\thetai_1, \ldots , \thetai_4$, with $i=1,2$, represent angles and, hence, are dimensionless by definition.
We further introduce the length-ratio $\alpha$, the respective slendernesses $\beta\ith{i}$ of the flexures, the width-ratio $\gamma$, and the horizontal position $X_0\ith{2}$ of the point $s=0$ of the second flexure, which is scaled by the length $l\ith{1}$ of the first flexure giving a non-dimensional ratio $\delta$: 
\begin{equation}
    \alpha = \frac{l\two}{l\one} , \qquad
    \beta\ith{i} = \frac{l\ith{i}}{h\ith{i}} , \qquad
    \gamma = \frac{w\ith{2}}{w\ith{1}} , \qquad
    \delta = \frac{X_0\ith{2}}{l\ith{1}} .
\end{equation}

To allow for effective optimization, we need to impose restrictions on the design variables that preclude rotated and mirrored shapes, which behave identical from in terms of objectives which are to be adopted subsequently.
For this purpose, we fix the point $\Xv\ith{2}(0) = \Xv_0\ith{2}$ of the second flexure to the $y$-axis if the global Cartesian frame, i.e., $Y_0\ith{2} = 0$. 
To avoid reflection symmetries, $\delta$ is required to be positive and the angle of the cross-section at $s=0$ of the first flexure is restricted to the range $\theta_0\ith{1} \in [0,\pi]$ 
The individual bounds for all \num{13} design variables as being used in the optimization are listed in Tab.~\ref{tab:range-variables}.
Note that infeasible shapes may occur when randomly sampling from the parameter ranges. 
In particular, self-intersecting centerlines of the flexures are not admissible and are rejected in the optimization process.

\section{PARETO-OPTIMIZATION OF CROSS-HINGE GEOMETRY} \label{sec:optimization}
In what follows, we assume the body attached to the flexures at $s=0$ as spatially fixed. 
The (rigid) body connected to cross-sections $s=1$ rotates by the angle $\varphi$ about the $z$-axis. 
The operational range (action space) of the cross-hinge is specified as $\varphi \in \mathcal A = \left[ 0, \pi / 2\right]$.
\subsection{Kinetostatic performance measures}
In what follows, we formulate the objective functions based on the above characterization of a hinge's nature.
Accordingly, we first introduce an appropriate measure to quantify the kinematic deviation of the cross-hinge from an ideal revolute joint.
Subsequently, static properties of the cross-hinge, i.e., its compliance with respect to translational motion and the stiffness against rotation are introduced.
%
\paragraph{Kinematics.} The distinguishing feature of a hinge's axis of rotation being stationary relative to the connected body is in the center of attention.
Within our approach, we refrain from specifying the actual position of the axis a priori.
Instead, we want to minimize the variation of the axis' position relative to connected bodies in the course of motion.
For this purpose, we resort to the kinematic concept of centrodes, specifically the \emph{fixed} centrode as the curve described by the instantaneous center of rotation relative to the fixed body.
Let $\xv\ith{i}(s, \varphi)$ denote the position of the centerlines for a prescribed (relative) rotation $\varphi$:
%
\begin{equation}
    \xv\ith{i}(s, \varphi) = \Xv\ith{i}(s) + \uv\ith{i}(s, \varphi) ,
    \qquad \varphi \in \mathcal A .
\end{equation}
We choose the terminal point $A$ of the flexure $i=1$ as a reference point on the rotating body, whose position as a function of the rotation $\varphi$ is given by: 
\begin{equation}
    \xv_A ( \varphi ) = \xv\ith{1}(s = 1, \varphi) .
\end{equation}
Recalling elementary kinematics of rigid bodies, the variation of the current position of $A$ can be expressed in terms of the instantaneous center of rotation $\xv_c$ as
\begin{equation}
    \frac{\partial \xv_A}{\partial \varphi} \delta \varphi = \delta \varphi \ez \times \left( \xv_A - \xv_c \right) .
\end{equation}
In the course of the optimation, we increase the prescribed angle $\varphi$ in discrete steps, which allows us to approximate the fixed centrode $\varphi \mapsto \xv_c$ replacing the derivative in the above relation by the corresponding difference quotient:
\begin{equation}
    \xv_c \approx \xv_A + \ez \times \frac{\Delta \xv_A}{\Delta \varphi} .
    \label{eq:centrode}
\end{equation}
While we used the arc-length of the centrode in previous studies~\cite{humer2018}, we propose to use the radius $r_c$ of the smallest circle containing the centrode instead.
Although both the length of the centrode and $r_c$ vanish ideally, the latter preserves information on the spatial variation of the axis as the cross-hinge deforms, whereas the former only captures the absolute distance by which the instant center of rotation get displaced. 
%
%
\paragraph{Translational compliance.} 
An ideal revolute joint is characterized by an infinite stiffness in all directions except for rotational degree of freedom. 
Even though conventional hinges hardly meet this requirement, the finite stiffness of compliant mechanisms against translational motion is typically significantly smaller.
When rotating the cross-hinge in discrete steps, we determine the principal translational compliances of the mechanism in each equilibrium configuration throughout the action space $\mathcal A$.
The largest of all principal compliances, i.e.,
\begin{equation}
    c_{\rm max} = \max_{j = 1,2 , \; \varphi \, \in \, \mathcal A} c_{{\rm transl}, j}(\varphi) ,
\end{equation}
serves as a measure to quantify the finite stiffness (or, much rather, the lack thereof) of the cross-hinge with respect to parasitic loads.

\paragraph{Rotational stiffness.}
Conversely, we seek for a minimum stiffness against a rotation about the indended axis.
In analogy to the translational motion, we consider the maximum rotational stiffness within the action space $\mathcal A$ as a measure which is to be minimized by the cross-hinge:
\begin{equation}
    k_{\rm min} = \max_{\varphi \, \in \, \mathcal A} k_{{\rm rot}} (\varphi) . 
\end{equation}

To obtain cross-hinge designs that can be easily scaled, we prefer to formulate the objectives as non-dimensional functions of the kinetostatic properties introduced above.
For this purpose, the length of the first flexure serves as a natural choice for the referential length, i.e., $\lref = l\ith{1}$.
To define a referential force, we take the product of Young's modulus with the referential length and the width of the first flexure, i.e., $\Fref = E l\ith{1} w\ith{1}$.
Accordingly, the set of objective functions is given by
\begin{equation}
    \rnd = \frac{r_c}{\lref}= \frac{r_c}{l\ith{1}} , \qquad 
    \cnd = \frac{c_{\rm max} \Fref}{\lref} = c_{\rm max} E w\ith{1}, \qquad
    \knd = \frac{k_{\rm min}}{\Fref \lref} = \frac{k_{\rm min}}{E w\ith{1} \left(l\ith{1}\right)^2 } ,
\end{equation}
where bars have been added to distinguish non-dimensional quantities. 
Referring stiffnesses and compliances to the width $w\ith{1}$, we essentially eliminate effects related to the depth of the mechanism $z$-direction from the objectives.
While this makes perfect sense when using (planar) beam theory in the analysis, such effects are certainly present in three-dimensional solid models of compliant mechanisms, as our subsequent results will show.

\subsection{Formulation of the optimization problem}
The individual performance measures are combined to form a vector-valued objective function:
\begin{equation}
    f : \D \rightarrow \Y \subseteq \reals^3 , \qquad x \mapsto \left( \rnd, \cnd, \knd \right)^T .
\end{equation}
In general, a global optimum that minimizes all components of $f$ does not exist.
From a mathematical point of view, we are lacking an order relation that allows us to compare two individuals $x$, $x^\prime$ in the sense $f(x) \leq f(x^\prime)$.
We adopt the concept of \emph{Pareto-optimality}, which has its roots in economy and, ever since, has successfully been applied to engineering problems~\cite{Deb2001,humer2016,humer2025}.

Let $x$, $x^\prime$ denote two solutions characterized by their respective objective vectors $y = f(x)$, $y^\prime = f(x^\prime)$.
Solution $x$ \emph{dominates} $x^\prime$ if $x$ is at least as good as $x^\prime$ with respect to all objectives and better than $x^\prime$ with respect to (at least) one objective, i.e.,
\begin{equation}
    \forall i \in \left\{ 1, \ldots, k \right\} : \; y_i \leq y_i^\prime  \quad \wedge \quad
    \exists	i \in \left\{ 1, \ldots, k \right\} : \; y_i < y_i^\prime .
\end{equation} 
To indicate the Pareto order, the short-hand notation
\begin{equation}
    y \prec y^\prime , \qquad x \prec_f x^\prime 
\end{equation}
is typically used.
The set of all \emph{non-dominated} solutions represents the \emph{Pareto set}, 
\begin{equation}
    X_P = \left\{ x \in \D : \left\{ x' \in \D : x' \prec_f x , x \neq x' \right\} = \emptyset \right\} ;
\end{equation}
the \emph{Pareto front} is the image of the Pareto set in the objective space, i.e.,
\begin{equation}
    Y_P = \left\{ y \in \Y : \left\{ y' \in \Y : y' \prec y , y \neq y' \right\} = \emptyset \right\} .
\end{equation}
Although the optimization problem is formulated in terms of non-dimensional objectives, their numerical values may differ by several orders of magnitude.
We normalize the Pareto front for the sake of comparability of the individual objectives.
To this end, we introduce the concepts of the the \emph{ideal} objective vector and the \emph{nadir} point, respectively, 
\begin{equation}
    \ideal = \left( \inf_{x' \in X_P} f_1(x') , \ldots , \inf_{x' \in X_P} f_k(x') \right)^T , \qquad
    \nadir = \left( \sup_{x' \in X_P} f_1(x') , \ldots , \sup_{x' \in X_P} f_k(x') \right)^T ,
    \label{eq:ideal_nadir}
\end{equation}
which represent the respective lower and upper bounds of the individual objectives over the Pareto set. 
The (approximations of the) ideal and nadir solutions are used to normalize the Pareto front, i.e.,
\begin{equation}
	\tilde Y_P = \left\{ \tilde y_i = \frac{y_i - \ideal}{\nadir - \ideal} : y_i \in Y_P \right\} \subseteq [0, 1]^3 ,
\end{equation}
where a tilde has been introduced to indicate normalized objectives subsequently.

\subsection{Implementation and solution}
The proposed hybrid structural-continuum approach relies two distinct models of the cross-hinge. 
To approximate the Pareto set by means of evoluationary algorithms, we employ the geometrically exact theory for the planar deformation of shear-deformable beams, which goes back to Reissner~\cite{reissner1972}.
We assume the conventional constitutive relations of structural mechanics, by which cross-sectional stress resultants and generalized strains are related linearly.
The absolute value of Young's modulus is not relevant within our non-dimensional framework.
Poisson's ratio is set to $\nu = \num{0.49}$ to account for nearly incompressible materials as used in additively manufacturing.
Although closed-form solution under concentrated loads do exist, see, e.g.,~\cite{humer2013, humer2019}, we use a finite-element model based on the formulation of Simo and Vu-Quoc~\cite{simo1986a}.
Each flexure of the cross-hinge is discretized into \num{30} cubic beam elements.
The ends at $s=0$ are fully constrained; cross-sections at $s=1$ share three degrees of freedom representing a rigid-body motion.
The rotation of $\pi/2$ at $s=1$ is kinematically imposed in \num{20} steps. 
After each step, we compute a Schur complement of the tangent stiffness matrix to determine the stiffnesses with respect to a translational motion of a body attached at $s=1$, from which the principal compliances follow as inverse of corresponding eigenvalues.
The rotational stiffness is obtained by a finite-difference approximation of the reaction moment that occurs in the course of deformation.
Welzl's algorithm~\cite{welzl1991} is used to compute the circumradius of the centrode~\eqref{eq:centrode}.
To distribute the compliance along the flexures, we impose a constraint on the bending strain, which is bound from above by \num{0.2}.
A single cross-hinge design can be analyzed within a few tenths of a second.
Both SPEA2~\cite{zitzler2001} and NSGA-II~\cite{deb2002} evolutionary algorithms as implemented in the open-source Python package \texttt{pymoo}~\cite{pymoo} are used to find Pareto-optimal solutions.

Once a sufficient approximation of the Pareto-set is obtained, selected solutions are to be refined using three-dimensional finite-element models of the cross-hinge. 
Specifically, we employ third-order elements, which almost directly replace the beam elements of the structural model.
A nearly incompressible, Neo-Hookean constitutive model governs the stress-response of the flexures.
The kinetostatic performance measures are evaluated similarly as for the structural model.
The computational effort to analyze a single cross-hinge design using the high-fidelity model exceeds the structural model by more than \num{1000} times, which is prohibitively large for evolutionary algorithms.
To fine-tune cross-hinge designs, we apply the Nelder-Mead~\cite{nelder1965} algorithm to a single-objective problem obtained by scalarization.
All finite-element models are implemented in the open-source code \href{https://www.ngsolve.org}{\texttt{Netgen/NGSolve}}.

\section{\MakeUppercase{Results and discussion}}
\subsection{Pareto-optimal solutions}
\begin{figure}[h!]
    \centering
    \includegraphics[scale=0.75]{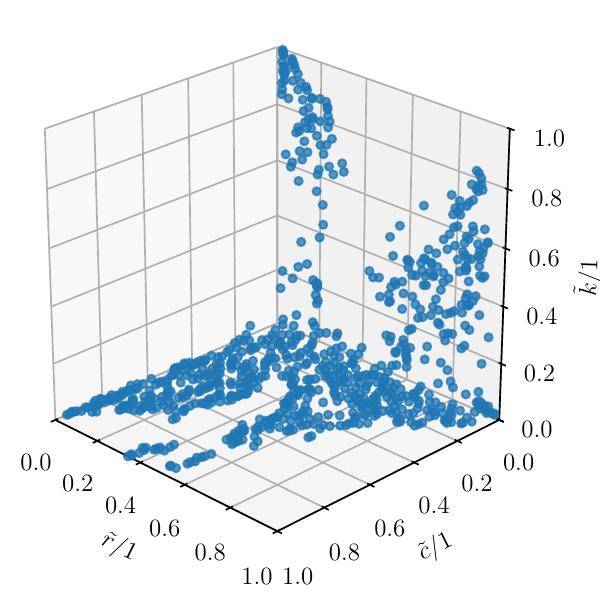}
    \caption{Approximation of the Pareto front comprising \num{798} individuals.}
    \label{fig:pareto_3d}
\end{figure}

To approximate the Pareto front based on the structural model, we evolve populations of \num{500} individuals over \num{1000} generations using SPEA2 and NSGA-II, respectively.
It turns out that the algorithms favor different directions of exploration for the present problem.
For this reason, we merge the respective Pareto sets, for which a total of \num{1000000} designs has been simulated. 
Figure~\ref{fig:pareto_3d} shows the approximation of Pareto front, which is a surface in the three-dimensional objective space. 
Despite the seemingly simple nature of the present problem, the Pareto front has a comparatively complex shape showing several disconnected regions, which already indicate the presence of diverse design variants. 

\begin{figure}[h!]
    \centering
    \includegraphics[scale=0.75]{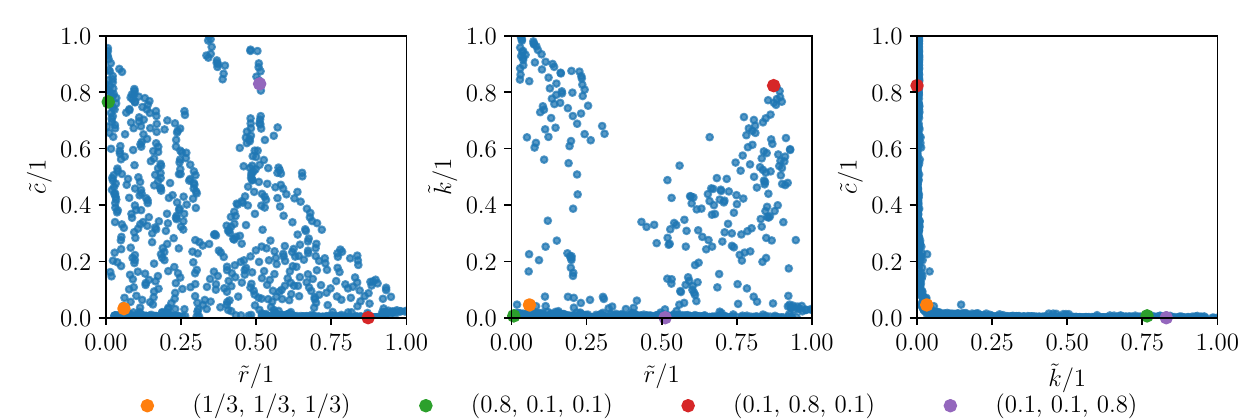}
    \vspace{-0.5cm}
    \caption{Projections of the Pareto front onto two-dimensional subspaces spanned by pairs of objectives. Colored markers correspond to selected solutions obtained upon uniform and preferential weighting, respectively.}
    \label{fig:pareto_2d}
\end{figure}

Projections of the Pareto front onto two-dimensional subspaces spanned by pairs of the objectives are illustrated in Fig.~\ref{fig:pareto_2d}.
As all objectives are to be minimized, solutions optimal with respect to the respective pairs of lie towards the lower left.
The projection onto the $\tilde c \tilde k$-plane (right) shows a typical trade-off among two conflicting objectives: A smaller rotational stiffness $k$ implies a greater translational compliance $c$ and vice versa.
Projections that include the deviation from a hinge's kinematics in terms of the centrode circumradius $r$, by contrast, show a more complex relationship among the objectives.
We find disconnected regions in projections both onto the $\tilde c \tilde r$-plane (left) and the $\tilde k \tilde r$-plane (center).

\subsection{Multi-criteria decision making}

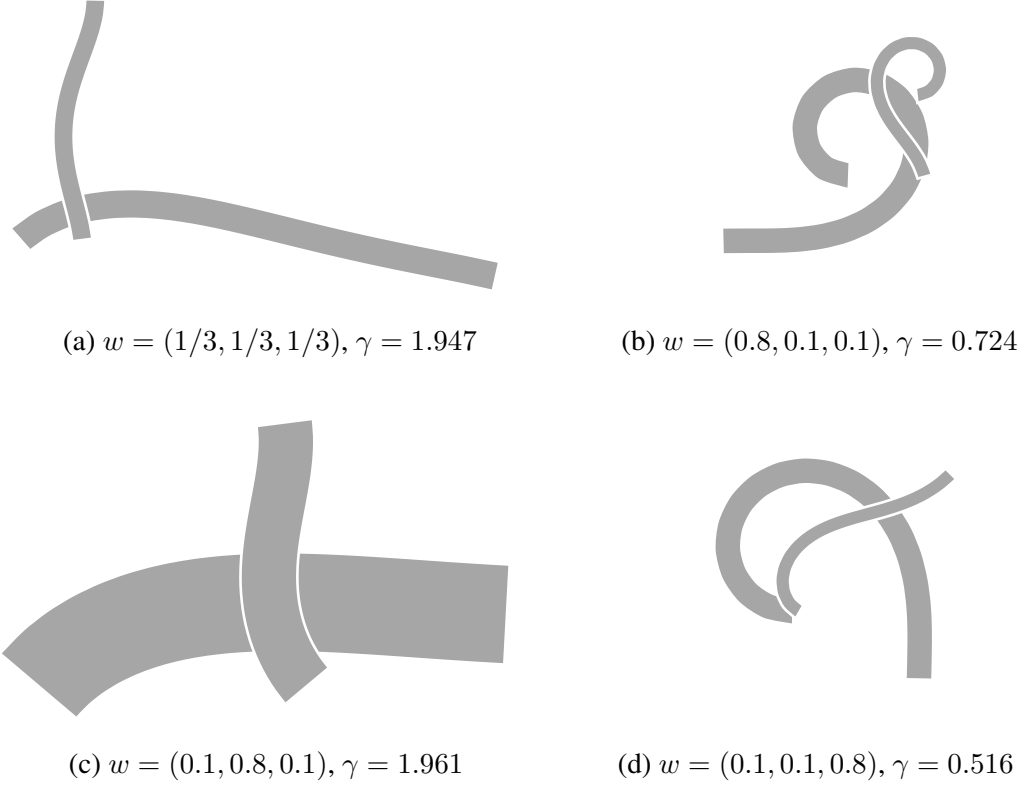
\begin{figure}[ht!]
    \centering
	\begin{minipage}[t][0.3\textwidth][c]{0.475\textwidth} 
		\centering
		\begin{tikzpicture}[scale=0.65]
            \drawSmoothThickCurve{(9.744,-0.4935) (9.25,-0.3867) (8.758,-0.287) (8.268,-0.1896) (7.779,-0.09108) (7.292,0.01111) (6.807,0.1183) (6.323,0.2307) (5.839,0.3475) (5.354,0.467) (4.867,0.5862) (4.376,0.7009) (3.879,0.806) (3.375,0.8948) (2.863,0.9591) (2.343,0.9894) (1.816,0.9745) (1.288,0.9018) (0.7677,0.758) (0.2697,0.5305) (-0.1836,0.2094) }{(0.1836,-0.2094) (0.5451,0.0463) (0.9557,0.2336) (1.398,0.3557) (1.86,0.4192) (2.333,0.4325) (2.811,0.4046) (3.291,0.3442) (3.772,0.2594) (4.253,0.1575) (4.736,0.04465) (5.22,-0.07376) (5.706,-0.1935) (6.194,-0.3113) (6.684,-0.425) (7.175,-0.5334) (7.667,-0.6367) (8.159,-0.7359) (8.649,-0.8333) (9.137,-0.9321) (9.621,-1.037) }
            \drawSmoothThickCurve{(1.334,4.875) (1.313,4.65) (1.265,4.423) (1.197,4.194) (1.117,3.963) (1.031,3.729) (0.9458,3.49) (0.8658,3.246) (0.7962,2.996) (0.7408,2.741) (0.7026,2.481) (0.6834,2.218) (0.6836,1.954) (0.7025,1.693) (0.7382,1.434) (0.7873,1.18) (0.8459,0.931) (0.9089,0.6869) (0.9705,0.447) (1.025,0.2105) (1.064,-0.02293) }{(1.423,0.02293) (1.379,0.2827) (1.322,0.5339) (1.259,0.7785) (1.197,1.019) (1.141,1.257) (1.095,1.494) (1.062,1.731) (1.045,1.968) (1.045,2.205) (1.062,2.441) (1.097,2.676) (1.147,2.909) (1.212,3.141) (1.288,3.372) (1.371,3.604) (1.457,3.84) (1.541,4.082) (1.615,4.332) (1.67,4.594) (1.695,4.868) }
        \end{tikzpicture}
	\end{minipage}
    \begin{minipage}[t][0.3\textwidth][c]{0.475\textwidth} 
        \centering
        \begin{tikzpicture}[scale=0.65]
            \drawSmoothThickCurve{(-2.539,-1.592) (-2.046,-1.588) (-1.541,-1.586) (-1.026,-1.564) (-0.5043,-1.497) (0.01474,-1.363) (0.5133,-1.142) (0.9642,-0.8192) (1.331,-0.3921) (1.572,0.1257) (1.645,0.6997) (1.522,1.272) (1.197,1.768) (0.7002,2.102) (0.1043,2.204) (-0.482,2.036) (-0.932,1.617) (-1.132,1.031) (-1.019,0.4181) (-0.6088,-0.05664) (-0.008751,-0.2505) }{(0.008751,0.2505) (-0.3306,0.3605) (-0.5661,0.6336) (-0.6322,0.9928) (-0.5123,1.343) (-0.2372,1.599) (0.1295,1.703) (0.5116,1.637) (0.8384,1.417) (1.058,1.081) (1.144,0.6812) (1.091,0.2676) (0.9113,-0.118) (0.6281,-0.4472) (0.2673,-0.705) (-0.147,-0.8886) (-0.5955,-1.004) (-1.065,-1.064) (-1.549,-1.085) (-2.044,-1.086) (-2.551,-1.09) }
            \drawSmoothThickCurve{(1.443,1.519) (1.733,1.628) (1.935,1.861) (2.006,2.16) (1.933,2.457) (1.738,2.69) (1.465,2.817) (1.167,2.824) (0.8908,2.718) (0.6728,2.524) (0.5312,2.273) (0.4696,1.996) (0.4804,1.717) (0.5497,1.452) (0.6612,1.207) (0.7993,0.9836) (0.9502,0.7764) (1.102,0.578) (1.242,0.3806) (1.359,0.1777) (1.441,-0.03334) }{(1.683,0.03334) (1.584,0.2883) (1.451,0.5182) (1.302,0.7286) (1.15,0.928) (1.006,1.125) (0.8811,1.327) (0.7851,1.537) (0.728,1.755) (0.7193,1.975) (0.7657,2.185) (0.8688,2.368) (1.022,2.505) (1.211,2.578) (1.409,2.573) (1.586,2.491) (1.71,2.343) (1.755,2.158) (1.712,1.975) (1.589,1.834) (1.415,1.768) }
        \end{tikzpicture}
    \end{minipage}
    \\
    \begin{minipage}[t][1.5em][t]{0.45\textwidth} \small \centering (a) $\weight = \left(1/3, 1/3, 1/3 \right)$, $\gamma = \num{1.947}$ \end{minipage}
	\begin{minipage}[t][1.5em][t]{0.45\textwidth} \small \centering (b) $\weight = \left(0.8, 0.1, 0.1 \right)$, $\gamma = \num{0.724}$ \end{minipage}
    \\
    \begin{minipage}[t][0.3\textwidth][c]{0.475\textwidth} 
        \centering
        \begin{tikzpicture}[scale=0.65]
            \drawSmoothThickCurve{(9.597,2.44) (9.105,2.469) (8.612,2.5) (8.116,2.534) (7.617,2.569) (7.115,2.603) (6.609,2.635) (6.097,2.662) (5.58,2.682) (5.055,2.691) (4.524,2.685) (3.986,2.66) (3.441,2.612) (2.891,2.534) (2.337,2.421) (1.782,2.267) (1.232,2.064) (0.6924,1.806) (0.1723,1.487) (-0.3171,1.101) (-0.7615,0.6455) }{(0.7615,-0.6455) (1.015,-0.3861) (1.31,-0.1533) (1.642,0.05005) (2.005,0.223) (2.392,0.3657) (2.801,0.4795) (3.227,0.5661) (3.666,0.6281) (4.117,0.6682) (4.578,0.6892) (5.047,0.6943) (5.523,0.6863) (6.004,0.668) (6.491,0.6421) (6.983,0.611) (7.478,0.5771) (7.977,0.5422) (8.479,0.508) (8.984,0.4758) (9.49,0.4467) }
            \drawSmoothThickCurve{(4.477,5.277) (4.493,5.057) (4.489,4.824) (4.468,4.581) (4.433,4.328) (4.389,4.068) (4.338,3.799) (4.287,3.521) (4.238,3.236) (4.196,2.943) (4.165,2.642) (4.148,2.334) (4.15,2.022) (4.174,1.706) (4.22,1.389) (4.293,1.075) (4.391,0.7657) (4.516,0.4652) (4.667,0.1762) (4.842,-0.09822) (5.04,-0.3558) }{(5.891,0.3558) (5.751,0.5387) (5.627,0.7321) (5.521,0.9353) (5.433,1.147) (5.363,1.367) (5.311,1.594) (5.276,1.826) (5.259,2.064) (5.257,2.307) (5.27,2.554) (5.296,2.805) (5.333,3.061) (5.378,3.323) (5.428,3.592) (5.48,3.869) (5.529,4.156) (5.57,4.453) (5.597,4.763) (5.602,5.085) (5.577,5.418) }
        \end{tikzpicture}
    \end{minipage}
    \begin{minipage}[t][0.3\textwidth][c]{0.475\textwidth} 
        \centering
        \begin{tikzpicture}[scale=0.65]
            \drawSmoothThickCurve{(2.846,-1.378) (2.857,-0.8776) (2.859,-0.3689) (2.835,0.1473) (2.766,0.6671) (2.638,1.183) (2.436,1.682) (2.15,2.145) (1.775,2.549) (1.318,2.865) (0.7935,3.065) (0.2304,3.124) (-0.3312,3.023) (-0.8416,2.759) (-1.248,2.347) (-1.5,1.823) (-1.562,1.242) (-1.416,0.6734) (-1.076,0.1923) (-0.5825,-0.1323) (-0.001852,-0.2501) }{(0.001852,0.2501) (-0.3915,0.3299) (-0.7283,0.5513) (-0.962,0.8823) (-1.063,1.277) (-1.019,1.685) (-0.8396,2.058) (-0.5463,2.355) (-0.1717,2.548) (0.2475,2.624) (0.676,2.579) (1.083,2.424) (1.446,2.172) (1.751,1.844) (1.989,1.458) (2.162,1.031) (2.275,0.5767) (2.336,0.1054) (2.359,-0.3771) (2.357,-0.8684) (2.346,-1.368) }
            \drawSmoothThickCurve{(3.139,2.882) (2.963,2.726) (2.77,2.589) (2.562,2.471) (2.344,2.371) (2.117,2.285) (1.885,2.21) (1.647,2.143) (1.406,2.077) (1.162,2.009) (0.9184,1.931) (0.6767,1.837) (0.4412,1.72) (0.2182,1.573) (0.01679,1.392) (-0.1505,1.174) (-0.2679,0.9209) (-0.3178,0.6416) (-0.2834,0.3552) (-0.1541,0.09195) (0.06729,-0.108) }{(0.1941,0.108) (0.047,0.2412) (-0.04276,0.4245) (-0.06753,0.6326) (-0.02958,0.844) (0.06295,1.043) (0.1996,1.221) (0.3699,1.374) (0.5647,1.502) (0.7769,1.607) (1.001,1.694) (1.233,1.769) (1.471,1.836) (1.713,1.901) (1.956,1.97) (2.199,2.048) (2.44,2.14) (2.676,2.248) (2.903,2.378) (3.118,2.53) (3.316,2.705) }
        \end{tikzpicture}
    \end{minipage}
    \\
    \begin{minipage}[t][1.5em][t]{0.45\textwidth} \small \centering (c) $\weight = \left(0.1, 0.8, 0.1 \right)$, $\gamma = \num{1.961}$ \end{minipage}
	\begin{minipage}[t][1.5em][t]{0.45\textwidth} \small \centering (d) $\weight = \left(0.1, 0.1, 0.8 \right)$, $\gamma = \num{0.516}$ \end{minipage}
    \caption{Schematic representations of cross-hinge geometries corresponding to uniform weighting among objectives (a) and preferred weighting of one of the objectives each (b)-(d).}
    \label{fig:shapes_2d}
\end{figure}
Having approximated the Pareto set, we are only halfway through in multi-objective optimization. 
We still face the question of how to systematically choose individual soluations from the Pareto set.
The field of \emph{multi-criteria decision making} provides us with different approaches, among which we employ the concept of \emph{pseudo-weights} as introduced by~\cite{Deb2001}.
For each solution of the Pareto set $x \in X_P$, a vector of pseudo-weights $\wpseudo \in \reals^3$ is computed based on the normalized distance to the worst solution, i.e., the nadir point $\nadir$, cf. Eq.~\eqref{eq:ideal_nadir}:
\begin{equation}
	\wpseudoi (x)
	= \frac{\nadiri - f_i(x)}{\nadiri - \ideali} \sum_{j=1}^3 \frac{\nadirj - \idealj}{\nadirj - f_j(x)}
	= \frac{1 - \tilde f_i(x)}{\sum_{j=1}^3 \left( 1 - \tilde f_j(x) \right)} , \qquad i = 1, \ldots , 3 .
	\label{eq:pseudo_weights}
\end{equation}
To make a selection, we choose the particular solution $x^*$ that minimizes the $L_1$-distance between a preferred (target) weighting $w = (w_r, w_c, w_d)^T$ and the vector of pseudo-weights $\wpseudo$:
\begin{equation}
	x^* = \argmin_{x \in X_P} \sum_{i=1}^3 \vert \wpseudoi (x) - w_i \vert .
	\label{eq:pseudo_weights-l1}
\end{equation}
Note that (pseudo) weights $w_i$ need to be non-negative and sum up to $1$, i.e., $\sum_{i=1}^3 w_i = 1$.

To provide some insight into the diversity within the set of Pareto-optimal cross-hinge designs, we consider four different target weightings:
A uniform weighting of $w_i = 1/3$, $i=1,2,3$ is meant to represent what could be referred to as ideal trade-off. 
Additionally, we consider cross-hinge designs that emphasize one of the objectives over the remaining two by setting $w_i = 0.8$ and $w_j = 0.1$ for all $j \neq i$, where $i,j \in \left\{ 1, 2, 3 \right\}$.
Figure~\ref{fig:shapes_2d} illustrates the shapes of the flexures corresponding to these weightings; the respective (normalized) objectives are listed in Tab.~\ref{tab:results_2d}.
The location of the specific designs are also highlighted by colored markers in the two-dimensional projections of the Pareto front, see Fig.~\ref{fig:pareto_2d}. 
\begin{table}[h!]
    \sisetup{
		round-mode=places,
		round-precision=3,
		output-exponent-marker = e
	}
    \centering
    \footnotesize
    \caption{Selected solutions of the Pareto set: uniform weighting of all objectives and preferred weighting of one of the objectives – target-weights, pseudo-weights and non-dimensional objectives.}
    \begin{tabular}{l|l|ccc|c}
        \toprule
        target weights  &   pseudo-weights  &   $\tilde r / 1$  &   $\tilde c / 1$   &  $\tilde k / 1$ &    shape \\
        \midrule
        (\num{0.333}, \num{0.333}, \num{0.333}) & 	 (\num{0.328}, \num{0.338}, \num{0.333}) & 	 \num{5.978e-02} & 	 \num{3.228e-02} & 	 \num{4.534e-02} & 	 (a) \\ 
        (\num{0.800}, \num{0.100}, \num{0.100}) & 	 (\num{0.447}, \num{0.105}, \num{0.447}) & 	 \num{7.513e-03} & 	 \num{7.658e-01} & 	 \num{6.467e-03} & 	 (b) \\ 
        (\num{0.100}, \num{0.800}, \num{0.100}) & 	 (\num{0.098}, \num{0.767}, \num{0.135}) & 	 \num{8.727e-01} & 	 \num{1.077e-04} & 	 \num{8.234e-01} & 	 (c) \\ 
        (\num{0.100}, \num{0.100}, \num{0.800}) & 	 (\num{0.295}, \num{0.103}, \num{0.603}) & 	 \num{5.115e-01} & 	 \num{8.298e-01} & 	 \num{6.610e-05} & 	 (d) \\
        \bottomrule
    \end{tabular}
    \label{tab:results_2d}
\end{table}

The cross-hinge designs show a remarkable variety in shape.
The shapes of the (undeformed) axes of the flexures obtained upon uniform weighting (a) show some similarity to the design preferring low translational compliance (c), but the latter design is characterized by much less slender flexures. 
The cross-hinge design preferring high kinematic accuracy (b) has two strongly curved, comparatively slender flexures. 
For a small rotational stiffness (d), the points of the two beams with respect to the ground body almost coincide.
All four designs share the property that flexure $\body_2$ has a smaller thickness than $\body_1$. 

The numerical results in Tab.~\ref{tab:results_2d} show that the individual performance measures vary by several orders of magnitude. 
It is also remarkable that target weights and pseudo-weights agree well for uniform weights (a) and small translational compliance (c). 
The mismatch between target and pseudeo-weights suggests that designs with pseudo-weights preferring high kinematic accuracy (b) or low rotational stiffness (d) in the desired extent are not present in the Pareto set.
 
\subsection{Fine-tuning based on three-dimensional model}
\begin{figure}[ht!]
    \centering
	\begin{minipage}[b][0.3\textwidth][c]{0.5\textwidth} 
		\centering
		\begin{tikzpicture}[scale=0.65]
            \drawSmoothThickCurveShadow{(9.744,-0.4935) (9.25,-0.3867) (8.758,-0.287) (8.268,-0.1896) (7.779,-0.09108) (7.292,0.01111) (6.807,0.1183) (6.323,0.2307) (5.839,0.3475) (5.354,0.467) (4.867,0.5862) (4.376,0.7009) (3.879,0.806) (3.375,0.8948) (2.863,0.9591) (2.343,0.9894) (1.816,0.9745) (1.288,0.9018) (0.7677,0.758) (0.2697,0.5305) (-0.1836,0.2094) }{(0.1836,-0.2094) (0.5451,0.0463) (0.9557,0.2336) (1.398,0.3557) (1.86,0.4192) (2.333,0.4325) (2.811,0.4046) (3.291,0.3442) (3.772,0.2594) (4.253,0.1575) (4.736,0.04465) (5.22,-0.07376) (5.706,-0.1935) (6.194,-0.3113) (6.684,-0.425) (7.175,-0.5334) (7.667,-0.6367) (8.159,-0.7359) (8.649,-0.8333) (9.137,-0.9321) (9.621,-1.037) }{gray!30}
            \drawSmoothThickCurveShadow{(1.334,4.875) (1.313,4.65) (1.265,4.423) (1.197,4.194) (1.117,3.963) (1.031,3.729) (0.9458,3.49) (0.8658,3.246) (0.7962,2.996) (0.7408,2.741) (0.7026,2.481) (0.6834,2.218) (0.6836,1.954) (0.7025,1.693) (0.7382,1.434) (0.7873,1.18) (0.8459,0.931) (0.9089,0.6869) (0.9705,0.447) (1.025,0.2105) (1.064,-0.02293) }{(1.423,0.02293) (1.379,0.2827) (1.322,0.5339) (1.259,0.7785) (1.197,1.019) (1.141,1.257) (1.095,1.494) (1.062,1.731) (1.045,1.968) (1.045,2.205) (1.062,2.441) (1.097,2.676) (1.147,2.909) (1.212,3.141) (1.288,3.372) (1.371,3.604) (1.457,3.84) (1.541,4.082) (1.615,4.332) (1.67,4.594) (1.695,4.868) }{gray!20}
            \drawSmoothThickCurve{(9.524,-1.584) (9.053,-1.407) (8.584,-1.234) (8.116,-1.061) (7.65,-0.8852) (7.186,-0.7058) (6.724,-0.5222) (6.262,-0.3348) (5.8,-0.1452) (5.335,0.04419) (4.867,0.23) (4.392,0.4078) (3.907,0.5718) (3.412,0.7148) (2.904,0.8282) (2.383,0.9017) (1.851,0.9233) (1.314,0.8795) (0.7823,0.7563) (0.272,0.5407) (-0.1918,0.2228) }{(0.1918,-0.2228) (0.5514,0.02332) (0.961,0.1962) (1.402,0.2981) (1.861,0.3354) (2.328,0.3164) (2.797,0.2501) (3.265,0.1456) (3.73,0.01119) (4.193,-0.1456) (4.654,-0.3182) (5.115,-0.501) (5.577,-0.6892) (6.04,-0.8791) (6.505,-1.068) (6.972,-1.253) (7.44,-1.435) (7.91,-1.612) (8.38,-1.786) (8.849,-1.959) (9.315,-2.133) }
            \drawSmoothThickCurve{(2.023,4.843) (1.935,4.623) (1.826,4.405) (1.703,4.19) (1.574,3.973) (1.444,3.753) (1.319,3.525) (1.205,3.29) (1.105,3.047) (1.023,2.795) (0.9613,2.537) (0.9218,2.275) (0.9042,2.01) (0.9072,1.746) (0.9282,1.485) (0.9633,1.227) (1.008,0.9736) (1.056,0.7242) (1.101,0.4785) (1.136,0.2357) (1.154,-0.003688) }{(1.446,0.003688) (1.427,0.2692) (1.39,0.5275) (1.343,0.7798) (1.295,1.028) (1.252,1.273) (1.219,1.517) (1.199,1.76) (1.197,2.002) (1.213,2.243) (1.249,2.481) (1.304,2.716) (1.379,2.945) (1.472,3.17) (1.579,3.39) (1.698,3.607) (1.825,3.823) (1.955,4.041) (2.083,4.266) (2.201,4.501) (2.3,4.75) }
        \end{tikzpicture}
	\end{minipage}
    \includegraphics[width=0.22\linewidth]{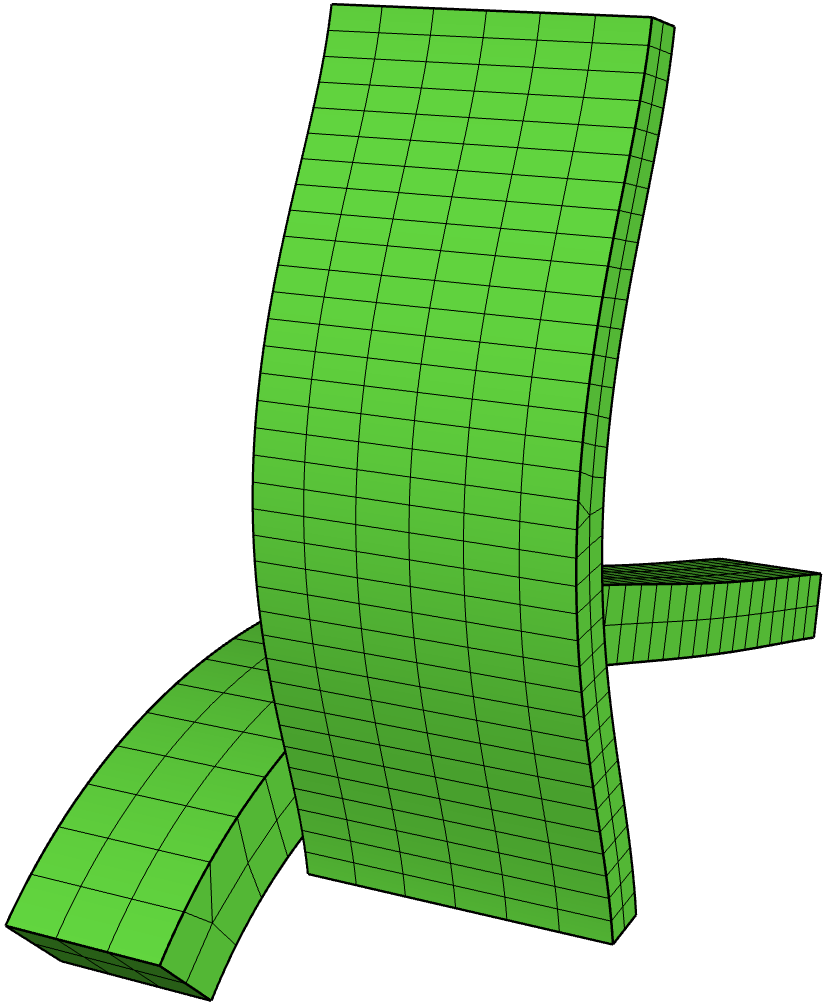}
    \caption{Schematic representation and three-dimensional finite-element model of the fine-tuned cross-hinge.}
    \label{fig:shape_3d}
\end{figure}

The structural model has proven a powerful tool to efficiently explore a the high-dimensional design space using evolutionary algorithms. 
Next, we want to fine-tune the cross-hinge designs using three-dimensional high-fidelity models. 
As an example, the shape obtained upon uniform weighting of the objectives, cf. Fig.~\ref{fig:shapes_2d} (a), is considered subsequently.
For this purpose, we scalarize the objective function using what is referred to as inverse normalization, i.e., the inverted normalized objectives are used as weights.
After \num{200} Nelder-Mead iterations, we obtain the cross-hinge design illustrated in Fig.~\ref{fig:shape_3d}, where the schematic representation is complemented by the three-dimensional finite-element model used in the analysis.  
For the sake of comparability, the original shape obtained from the structural model, which serves as initial solution for the scalar optimization, is included in light gray. 
Clearly, the original shape of the cross-hinge is mostly preserved; the adaptations through the single objective optimization using the three-dimensional model indeed improve the cross-hinge design with respect to all three objectives, see Tab.~\ref{tab:results_3d}.
The scalarized objective function $\tilde f_{\rm scalar} / 1$, which implies a weighting of $w_{\rm SOO} = (0.240, 0.444, 0.316)^T$, is reduced by almost \qty{20}{\percent} in the course of the single-objective optimization based on the high-fidelity model. 

\begin{table}[h!]
    \sisetup{
		round-mode=places,
		round-precision=3,
		output-exponent-marker = e
	}
    \centering
    \footnotesize
    \caption{Comparison of the structural model and the fine-tuned three-dimensional model of the cross-hinge: individual objectives and scalarized objective function.}
    \begin{tabular}{l|ccc|c}
        \toprule
                            &   $\tilde r / 1$  &   $\tilde c / 1$   &  $\tilde k / 1$ &    $\tilde f_{\rm scalar} / 1$ \\
        \midrule
        structural model    & 	 \num{5.978e-02} & 	 \num{3.228e-02} & 	 \num{4.534e-02} & 	 \num{4.300e-02} \\ 
        high-fidelity model & 	 \num{3.747e-02} & 	 \num{3.072e-02} & 	 \num{3.948e-02} & 	 \num{3.511e-02} \\ 
        \bottomrule
    \end{tabular}
    \label{tab:results_3d}
\end{table}

\section{\MakeUppercase{Conclusion}}
In the present paper, we have addressed the synthesis of compliant mechanisms with distributed compliances for applications characterized by large deformations.
A problem as simple as a hinge already provides some intution on the manifold aspect that need to be accounted for, not least from a computational point of view. 
Considering a generalized cross-hinge described by a set of \num{13} design variables, we have shown that the framework of multi-objective optimization provides effective means to explore feasible designs. 
Allowing for trade-offs among the kinetostatic objectives adopting the notion of Pareto-optimality, we have found are great variety of different shapes of the cross-hinge, which is reflected in the complex shape of the Pareto front. 
We have outlined the concept of pseudo-weights to systematically select individuals from the set of Pareto-optimal solutions.
Clearly, computational cost is key when using evolutionary algorithms.
For this reason, we have adopted a two-stage approach: Structural models based on the geometrically exact theory for shear-deformable beams have been used in the multi-objective optimization, where hundreds of thousands of variants are simulated efficiently. 
Subsequently, high-fidelity, high-order finite-element have been employed to further refine and fine-tune cross-hinge designs using  scalarized objective functions suited for conventional optimization methods.
\bibliographystyle{IEEEtranN}
\bibliography{References}
\end{document}